\documentclass[11pt, a4paper]{article}
\usepackage{array}
\usepackage{multicol}
\usepackage{multirow}
\usepackage{graphics}
\usepackage{graphicx}
\usepackage{epsfig}
\usepackage{amssymb}
\usepackage{amsthm}
\usepackage{mathrsfs}
\usepackage{amsmath}
\usepackage{float}
\usepackage{cite}
\usepackage{enumerate}
\usepackage{xcolor}
\usepackage{tabu}
\makeatletter
\newcommand{\thickhline}{%
    \noalign {\ifnum 0=`}\fi \hrule height 1pt
    \futurelet \reserved@a \@xhline
}
\newcolumntype{"}{@{\hskip\tabcolsep\vrule width 1pt\hskip\tabcolsep}}
\makeatother

\def\ms{\medskip}
\def\nt{\noindent}

\def\cyan{\color{cyan}}

\definecolor{vividviolet}{rgb}{0.62, 0.0, 1.0}
\def\rsq{\hspace*{\fill}$\Box$\medskip}

\newtheoremstyle{de}
  {10pt}          
  {10pt}  
  {\rm}  
  {}
  {\bf}  
  {. }    
  { }    
  {}     
\theoremstyle{de}

\newtheorem{example}{Example}[section]

\newtheorem{problem}{Problem}[section]

\newtheoremstyle{theorem}
  {10pt}          
  {10pt}  
  {\it}  
  {}
  {\bf}  
  {. }    
  { }    
  {}     
\theoremstyle{theorem}
\newtheorem{theorem}{Theorem}[section]
\newtheorem{lemma}[theorem]{Lemma}

\newtheorem{conjecture}{Conjecture}[section]

\numberwithin{equation}{section}

\graphicspath{%
    {converted_graphics/}
    {/}
}

\begingroup\makeatletter\ifx\SetFigFont\undefined%
\gdef\SetFigFont#1#2#3#4#5{%
  \reset@font\fontsize{#1}{#2pt}%
  \fontfamily{#3}\fontseries{#4}\fontshape{#5}%
  \selectfont}%
\fi\endgroup%

\begin{document}
\baselineskip18truept
\normalsize
\begin{center}
{\mathversion{bold}\Large \bf Construction of local antimagic 3-colorable graphs \\of fixed even size | matrix approach}

\bigskip
{\large Gee-Choon Lau$^{a,}$\footnote{Corresponding author.}, Wai Chee Shiu$^b$, M. Nalliah$^c$, K. Premalatha$^d$}

\medskip

\emph{{$^a$}77D, Jalan Subuh,\\
85000 Johor, Malaysia.}\\
\emph{geeclau@yahoo.com}

\medskip
\emph{{$^b$}Department of Mathematics,}\\
\emph{The Chinese University of Hong Kong,}\\

\emph{Shatin, Hong Kong, P.R. China.}\\
\emph{wcshiu@associate.hkbu.edu.hk}

\medskip
\emph{{$^c$}Department of Mathematics, School of Advanced Sciences, }\\
\emph{Vellore Institute of Technology, Vellore,\\ 632014, Tamil Nadu, India.}\\
\emph{nalliah.moviri@vit.ac.in}

\medskip
\emph{{$^d$}Department of Mathematics,}\\
\emph{Sri Shakthi Institute of Engineering and Technology,}\\
\emph{Coimbatore, 641062, Tamil Nadu, India.}\\
\emph{premalatha.sep26@gmail.com}
\end{center}

\begin{abstract}
An edge labeling of a connected graph $G = (V, E)$ is said to be local antimagic if it is a bijection $f:E \to\{1,\ldots ,|E|\}$ such that for any pair of adjacent vertices $x$ and $y$, $f^+(x)\not= f^+(y)$, where the induced vertex label $f^+(x)= \sum f(e)$, with $e$ ranging over all the edges incident to $x$.  The local antimagic chromatic number of $G$, denoted by $\chi_{la}(G)$, is the minimum number of distinct induced vertex labels over all local antimagic labelings of $G$. Suppose $\chi_{la}(G)=\chi_{la}(H)$ and $G_H$ is obtained from $G$ and $H$ by merging some vertices of $G$ with some vertices of $H$ bijectively. In this paper, we give ways to construct matrices with integers in $[1,10k]$, $k\ge 1$, that meet certain properties. Consequently, we obtained many families of (disconnected) bipartite (and tripartite) graphs of size $10k$ with local antimagic chromatic number 3.\\


\end{abstract}

\baselineskip18truept
\normalsize

\section{Introduction}
Let $G=(V, E)$ be a connected graph of order $n$ and size $m$.
A bijection $f: E\rightarrow \{1, 2, \dots, m\}$ is called a \textit{local antimagic labeling}
if $f^{+}(u)\neq f^{+}(v)$ whenever $uv\in E$,
where $f^{+}(u)=\sum_{e\in E(u)}f(e)$ and $E(u)$ is the set of edges incident to $u$.
The mapping $f^{+}$ which is also denoted by $f^+_G$ is called a \textit{vertex labeling of $G$ induced by $f$}, and the labels assigned to vertices are called \textit{induced colors} under $f$.
The \textit{color number} of a local antimagic labeling $f$ is the number of distinct induced colors under $f$, denoted by $c(f)$.  Moreover, $f$ is called a \textit{local antimagic $c(f)$-coloring} and $G$ is local antimagic $c(f)$-colorable. The \textit{local antimagic chromatic number} $\chi_{la}(G)$ is defined to be the minimum number of colors taken over all colorings of $G$ induced by local antimatic labelings of $G$. Let $G+H$ and $mG$ denote the disjoint union of graphs $G$ and $H$, and $m$ copies of $G$, respectively. For integers $a < b$, let $[a,b] = \{a, a+1, \ldots, b\}$.

\ms\nt The use of matrices that satisfy various properties have been used to determine the local antimagic chromatic number of many standard graphs (see~\cite{LNS, LPAS, LS, LSN}). However, the matrix used either has empty cells or repeated entries. It is interesting to obtain matrices without empty cells nor repeated entries such that the entries correspond to a local antimagic $\chi_{la}(G)$-coloring of a graph $G$. In~\cite{LSPN}, the authors used matrices of size $5\times (2k+1)$ and $11\times (2k+1)$ to construct various families of (disconnected) tripartite graphs of odd size $5\times (2k+1)$ and $11\times (2k+1)$  with local antimagic chromatic number 3.  In this paper, we give ways to construct  various matrices with even number of entries that meet certain properties.  Consequently, we obtained many new (disconnected) graphs (both bipartite and tripartite) of even size with local antimagic chromatic number 3. 

\ms\nt We shall need the following lemma. 

\begin{lemma}[\hspace*{-1mm}\cite{LSN}]\label{lem-2part}
Let $G$ be a graph of size $q$. Suppose there is a local antimagic labeling of $G$ inducing a $2$-coloring of $G$ with colors $x$ and $y$, where $x<y$.  Let $X$ and $Y$ be the sets of vertices colored $x$ and $y$, respectively, then $G$ is a bipartite graph with bipartition $(X,Y)$ and $|X|>|Y|$. Moreover,
\begin{equation*} x|X|=y|Y|= \frac{q(q+1)}{2}.\end{equation*}
\end{lemma}

\section{$5\times 2k$ matrix}

In this section, we shall always refer to the following $5\times 2k$ matrix (with entries in $[1,10k]$ bijectively) to get our results.

\[\fontsize{8}{11}\selectfont
\begin{tabu}{|c|[1pt]c|c|c|c|c|c|[1pt]c|c|c|c|c|c|}\hline
i & 1 & 2 & 3 & \cdots & k-1 & k & k+1 & k+2 & \cdots & 2k-2 & 2k-1 & 2k \\\tabucline[1pt]{-}
u_iw_i & 1 & k+1 & k+2 & \cdots & 2k-2 & 2k-1 & 2 & 3 & \cdots & k-1 & k & 2k \\\hline
v_iw_i & 6k & 6k+1 & 6k+2 & \cdots & 7k-2 & 7k-1 & 7k+1 & 7k+2 & \cdots & 8k-2 & 8k-1 & 8k \\\hline
x_iw_i & 7k & 6k-1 & 6k-3 & \cdots & 4k+5 & 4k+3 & 6k-2 & 6k-4 & \cdots & 4k+4 & 4k+2 & 3k+1 \\\hline
x_iu_i & 10k & 9k & 9k-1 & \cdots & 8k+3 & 8k+2 & 10k-1 & 10k-2 & \cdots & 9k+2 & 9k+1 & 8k+1 \\\hline
x_iv_i & 4k+1 & 4k & 4k-1 & \cdots & 3k+3 & 3k+2 & 3k & 3k-1 & \cdots & 2k+3 & 2k+2 & 2k+1 \\\hline
\end{tabu}\]

Observe that
\begin{enumerate}[(1)]
  \item In each column, the sum of the first 3 row entries is a constant $K_1=13k+1$.
  \item In each column, the sum of the rows 1 and 4 (respectively, 2 and 5) entries is $K_2=10k+1$. 
  \item The sums of the last 3 row entries for columns 1 and $2k$ are $21k+1$ and $13k+3$ respectively, while the corresponding sum for columns 2 to $k$ (respectively, $k+1$ to $2k-1$) form an arithmetic progression from $19k-1$ to $15k+7$ (respectively, from $19k-3$ to $15k+5$) in decrement of 4. Moreover, the corresponding sum for columns $a$ and $2k+1-a$ is $34k+4$ $(1\le a\le k)$.
  \item The sum of all the last 3 row entries is $S = k(34k+4)$.
  \item If $2k=rs, r \ge 2, s\ge 1$ then we can divide the table into $r$ blocks of $s$ column(s) with the $j$-th block containing $(j-1)s+1, (j-1)s+2, \ldots,js$ columns. Suppose $r$ is even. The sum of the row 3 entries in the $j$-th block, and the rows 4 and 5 entries in the $(r+1-j)$-th block $(1\le j\le r/2)$ is $K_3=s(17k+2)$. Similarly when $r$ is odd and $1\le j\le (r-1)/2$. Moreover, the sum of the last three row entries in $(r+1)/2$ block is also $K_3$.
  \item For $1\le i\le 2k$, the sum of the rows 2 and 3 entries of column $i$ and the row 4 entry of column $2k+1-i$ is a constant $K_4=21k+1$. 
  \item For $1\le i\le k$, the sum of row 4 columns $i$ and $2k+1-i$ entries is a constant $18k+1$, while the sum of row 5 columns $i$ and $2k+1-i$ entries is a constant $6k+2$.
\end{enumerate}

\nt Denote $nP_3\vee K_1$ by $FB(n)$, the {\it fan graph with $n$ blades}. Note that $FB(1)=P_3\vee K_1 \cong K_{1,2}\vee K_1$ is also the fan graph $F_3$ of order 4. In~\cite[Theorem 2.1]{LSPN}, the authors have proved that $]\chi_{la}(FB(n)) = 3$ for odd $n\ge 1$. We now extend it to even $n\ge 2$.

\begin{theorem}\label{thm-FBneven} For $k\ge 1$, $\chi_{la}(FB(2k)) = 3$.\end{theorem}

\begin{proof} Note that $FB(2k)$ is of size $10k$. Let the vertex set and edge set of the $i$-th copy of $FB(1)$ be $\{u_i, v_i, w_i, x_i\}$ and $\{u_iw_i, v_iw_i, x_iw_i, x_iu_i, x_iv_i\}$, respectively, $1\le i\le 2k$. Observe that table above gives a bijective edge labeling $f$ of $2k$ copies of $FB(1)$ using integers in $[1,10k]$ with induced vertex labels $f^+(u_i) = f^+(v_i) = 10k+1$, $f^+(w_i) = 13k+1$ and $\sum\limits^{2k}_{i=1} f^+(x_i) = k(34k+4)$ for $1\le i\le 2k$. Merging the vertices $x_1$ to $x_{2k}$ gives us $FB(2k)$ that has a vertex $x$ with $f^+(x)=k(34k+4)$. Thus, $\chi_{la}(FB(2k))\le 3$. Since $\chi_{la}(FB(2k))\ge \chi(FB(2k))=3$, the theorem holds.
\end{proof}

\begin{example}\label{eg-FB12} Take $k=6$, the following table
\[\fontsize{8}{11}\selectfont
\begin{tabu}{|c|[1pt]c|c|c|c|c|c|[1pt]c|c|c|c|c|c|}\hline
i & 1 & 2 & 3 & 4 & 5 & 6 & 7 & 8 & 9 & 10 & 11 & 12\\\tabucline[1pt]{-}
u_iw_i & 1 & 7 & 8 & 9 & 10 & 11 & 2 & 3 &  4 & 5  &6 & 12 \\\hline
v_iw_i & 36 & 37 & 38 & 39 & 40 & 41 & 43 & 44 & 45 & 46 & 47 & 48 \\\hline
x_iw_i & 42 & 35 & 33 & 31 & 29 & 27 & 34 & 32 & 30 & 28 & 26 & 19 \\\hline
x_iu_i & 60 & 54 & 53 & 52 & 51 & 50 & 59 & 58 & 57 & 56 & 55 & 49 \\\hline
x_iv_i & 25 & 24 & 23 & 22 & 21 & 20 & 18 & 17 & 16 & 15 & 14 & 13 \\\hline
\end{tabu}\]
we have a local antimagic labeling for $12FB(1)$ shown as follows:
\begin{figure}[H]
\centerline{\epsfig{file=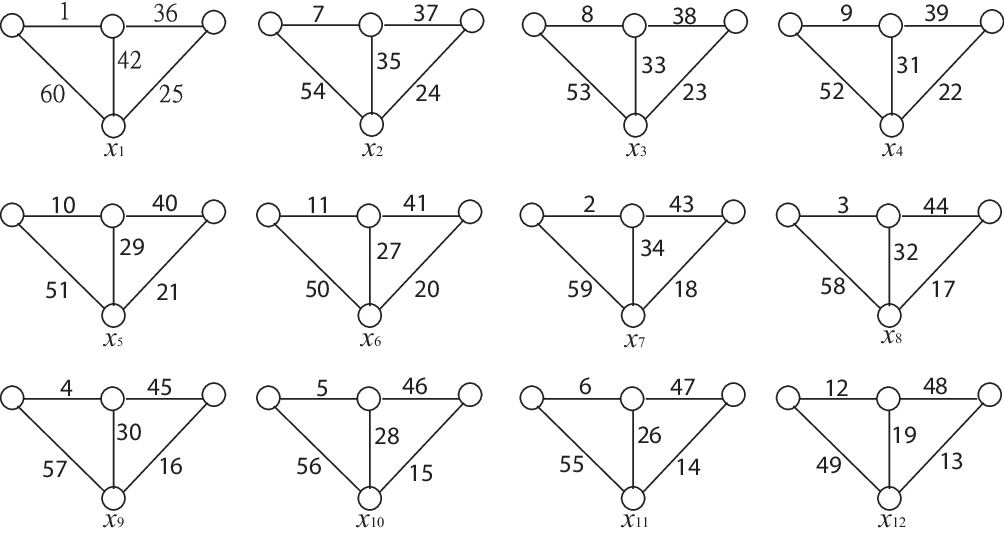, width=8.5cm}}
\caption{Twelve graphs that correspond to the 12 columns above.}\label{fig:FB12}
\end{figure}
\nt Merging vertices $x_1$ to $x_{12}$ gives the vertex $x$ of $FB(12)$ and the corresponding local antimagic $3$-coloring. \rsq
\end{example}

\begin{theorem}\label{thm-rFBseven} For $ rs\ge 4$ with $r\ge 2$ and even $s\ge 2$, $\chi_{la}(rFB(s)) = 3$. \end{theorem}

\begin{proof} Begin with $2k=rs$ copies of $FB(1)$ with edge labeling as in the proof of Theorem~\ref{thm-FBneven} and keep all the notation.

\nt For each $j \in[1,r]$, we merge the vertices in $\{x_{(j-1)s/2+a}, x_{2k-(j-1)s/2+1-a} \;|\; 1\le a\le s/2\}$.  We immediately get $r$ copies of $FB(s)$ with a bijective edge labeling using integers in $[1,10k]$. By Observation (3) above, $rFB(s)$ admits a local antimagic 3-coloring with
induced vertex labels $10k+1$, $13k+1$, $s(17k+2)$. Thus, $\chi_{la}(r(FB(s))) \le 3$. Since $\chi_{la}(r(FB(s)))\ge \chi(r(FB(s)))=3$, the theorem holds. \end{proof}

\begin{example} Consider $n=12$, we can obtain $3FB(4)$ (respectively, $6FB(2)$ and $2FB(6)$) from $12FB(1)$ with corresponding edge labeling similar to that in Example~\ref{eg-FB12} by merging the vertices in $\{x_{2i-1},x_{2i},x_{13-2i},x_{14-2i}\}$ for $1\le i\le 3$ (respectively, in $\{x_i,x_{13-i}\}$ for $1\le i\le 6$, and in $\{x_{3i-2},x_{3i-1},x_{3i},x_{13-3i},x_{14-3i},x_{15-3i}\}$ for $1\le i\le 2$). \rsq
\end{example}

\ms\nt For $r,s\ge 2$, $1\le i\le r$, let $G_i\cong FB(s)$ with vertex set $\{u_{i,j}, v_{i,j}, w_{i,j}, x_i\;|\;1\le j\le s\}$ and edge set $\{u_{i,j}w_{i,j}, v_{i,j}w_{i,j}, x_iw_{i,j}\;|\;1\le j\le s\}$. Note that $u_{i,j}$ and $v_{i,j}$ are of degree 2, $w_{i,j}$ is of degree 3, and $x_i$ is of degree $3s$.

\ms\nt Denote by $FB_1(r,s)$ the graph obtained from $G_1, \dots, G_r$ by merging the vertices in $\{u_{i,j}\;|\; 1\le i\le r\}$ and in $\{v_{i,j}\;|\; 1\le i\le r\}$ respectively for each $1\le j\le s$. Let the new vertices be $u_j$ and $v_j$ respectively for $1\le j\le s$. Note that $FB_1(r,s)$ has $rs$ vertices of degree 3, $2s$ vertices of degree $2r$, and $r$ vertices of degree $3s$.

\ms\nt Denote by $FB_2(r,s)$ graph obtained from $G_1, \dots, G_r$ by merging the vertices in $\{w_{i,j}\;|\; 1\le i\le r\}$ respectively for $1\le j\le s$. Let the new vertices be $w_j$ for $1\le j\le s$. Note that $FB_2(r,s)$ has $2rs$ vertices of degree 2, $s$ vertices of degree $3r$, and $r$ vertices of degree $3s$.

\begin{theorem}\label{thm-FBrs12} For $r\ge 2$ and even $s\ge 2$, $\chi_{la}(FB_1(r,s)) =3$ if $r\not\equiv 0\pmod 4$ and  $\chi_{la}(FB_2(r,s))=3$ if $rs\not\equiv 0\pmod 4$. \end{theorem}

\begin{proof} Let $rs = 2k$ for even $s\ge 2$. Since $\chi_{la}(FB_l(r,s))\ge \chi(FB_l(r,s))=3$ for $l=1,2$, we suffice to show that $\chi_{la}(FB_l(r,s))\le 3$. Using the $r$ copies of $FB(s)$ in Theorem~\ref{thm-rFBseven} and the corresponding local antimagic 3-coloring by merging the degree 2 vertices as defined for $FB_1(r,s)$, we can immediately conclude that $FB_1(r,s)$ admits a local antimagic 3-coloring with each $w_{i,j}$ (respectively, $u_j$, $v_j$ and $x_{i}$) has induced vertex label $13k+1$ (respectively, $r(10k+1)$ and $s(17k+2)$), for $1\le j\le s$ and $1\le i\le r$.

\ms\nt Suppose $r(10k+1)=s(17k+2)$.  If $r$ is odd, then we have $r^2(10k+1)=rs(17k+2)=2k(17k+2)$. So, $2k+1\equiv 2k^2\pmod 4$ which is impossible. If $r$ is even, then $k$ is even. Since $s$ is even, $r(10k+1)=s(17k+2)$ implies that $r\equiv 2s\equiv 0\pmod 4$. Thus if $r\not\equiv 0\pmod 4$, then $r(10k+1)\ne s(17k+2)$.

\ms\nt Similarly, by merging the degree 3 vertices as defined for $FB_2(r,s)$, we can immediately conclude that $FB_2(r,s)$ admits a local antimagic 3-coloring with each $u_{i,j}$, $v_{i,j}$ (respectively, $w_j$ and $x_i$) has induced vertex label $10k+1$ (respectively, $r(13k+1)$ and $s(17k+2)$), for $1\le j\le s$ and $1\le i\le r$.

\nt Suppose $r(13k+1)=s(17k+2)$. This implies that $2k(13k+1)=s^2(17k+2)$. Then $2k\equiv 0\pmod 4$. Thus when  $rs=2k\not\equiv 0\pmod 4$, $r(13k+1)\ne s(17k+2)$.

\ms\nt This completes the proof.
\end{proof}

\nt For $s\ge 1$, let $S_1$ and $S_2$ be two copies of $sP_3$. Let $DF(2s)$ be the {\it diamond fan} graph obtained from $S_1+S_2$ by joining a vertex $y$ (resp., $z$) to every degree 1 vertex of $S_1$ (resp., $S_2$) and every degree 2 vertex of $S_2$ (resp., $S_1$).  Thus, $DF(2s)$ has $4s$ vertices of degree 2, $2s$ vertices of degree 3, and 2 vertices of degree $3s$ with size $10s$.

\nt For $r\ge 1$, let $DF_r(2s) = rDF(2s) + FB(s)$. Note that $DF_r(2s)$ is of size $5(2r+1)s$ has
\begin{enumerate}[1.]
\item $(4r+2)s$ vertices of degree 2; these vertices will be denoted by $u_i$ and $v_i$, $1\le i\le (2r+1)s$;
\item $(2r + 1)s$ vertices of degree 3; these vertices will be denoted by $w_i$, $1\le i\le (2r+1)s$;
\item $2r+1$ vertices of degree $3s$; these vertices will be denoted by $x$,  $y_i$ and $z_i$, $1\le i\le r$. Hence $s$ can be 1.
\end{enumerate}

\begin{theorem}\label{thm-rDF2seven} For $r, s\ge 1$,  $\chi_{la}(rDF(2s))=3$.  \end{theorem}

\begin{proof} Let $2rs = 2k\ge 4$.  We first note that $rDF(2s)$ is a bipartite graph with equal partite set size. By Lemma~\ref{lem-2part}, $\chi_{la}(rDF(2s))\ge 3$. Suffice to show that $rDF(2s)$ admits a local antimagic 3-coloring.  

\ms\nt Begin with $2rs = 2k$ copies of $FB(1)$ with edge labeling as in the proof of Theorem~\ref{thm-FBneven}. Partition the $(2rs)FB(1)$ into $2r$ blocks of $s\ge 1$ copies of $FB(1)$ such that the $j$-th block has vertices $x_{(j-1)s+a}$ for $1\le j\le 2r, 1\le a\le s$. Split each $x_i$ into $x_i^1$ and $x_i^2$ such that $x_i^1$ is adjacent to $w_i$ and $x_i^2$ is adjacent to $u_i,v_i$. For each $1\le j\le r$ and $1\le a\le s$, merge the vertices in $\{x_{(j-1)s+a}^1, x_{(2r-j)s+a}^2\}$ to get vertex $y_j$, and the vertices in $\{x_{(j-1)s+a}^2, x_{(2r-j)s+a}^1\}$ to get vertex $z_j$. We now have $r$ copies of $DF(2s)$. By Observation (5),  we conclude that $rDF(2s)$ admits a local antimagic 3-coloring with degree  2 (respectively, $3$ and $3s$) vertices having induced vertex labels $10k+1$ (respectively $13k+1$ and $s(17k+2)$).  Thus, $\chi_{la}(rDF(2s)) \le 3$. This completes the proof.
\end{proof}

\begin{example}\label{eg-3DF4} Using the twelve copies of $FB(1)$ as in Example~\ref{eg-FB12} so that $k=6$, we can take $r=3, s=2$ to get $3DF(4)$ as shown below.
\begin{figure}[H]
\centerline{\epsfig{file=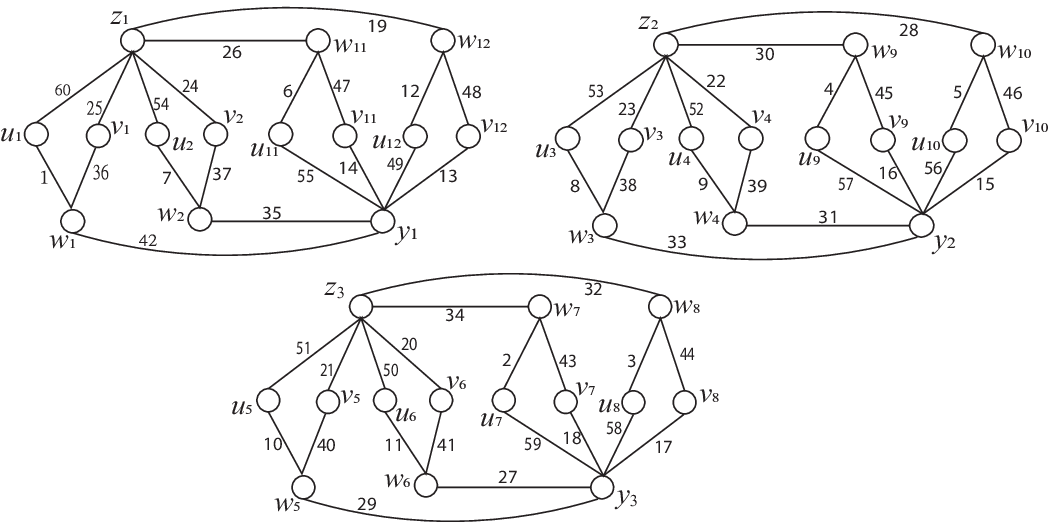, width=11.2cm}}
\caption{The $3DF(4)$ with the defined edge labeling.}\label{fig:3DF4}
\end{figure}
\nt Every degree 2 vertex has induced label 61, degree 3 vertex has induced label 79, and degree 6 vertex has induced label 208.
\rsq
\end{example}

\begin{theorem}\label{thm-DFr2seven} For $r\ge 1$ and even $s\ge 2$, $\chi_{la}(DF_r(2s))=3$.  \end{theorem}

\begin{proof} Let $(2r+1)s = 2k\ge 6$. Now, $DF_r(2s)$ has size $10k$. Since the $FB(s)$ component is a tripartite graph, we have $\chi_{la}(DF_r(2s))\ge \chi(DF_r(2s))=3$. Suffice to show that $DF_r(2s)$ admits a local antimagic 3-coloring.

\ms\nt Begin with $(2r+1)s$ copies of $FB(1)$ with edge labeling as in the proof of Theorem~\ref{thm-FBneven}. Partition the $(2r+1)sFB(1)$ into $2r+1$ blocks of even $s\ge 2$ copies of $FB(1)$ such that the $j$-th block has vertices $x_{(j-1)s+a}$ for $1\le j\le 2r+1, 1\le a\le s$. By the similar approach as in the proof of Theorem~\ref{thm-rDF2seven}, we can construct the $rDF(2s)$ using the $j$-th blocks for $j\in[1,2r+1]\setminus\{r+1\}$. The $(r+1)$-st block is used to construct the $FB(s)$ by merging vertices $x_{rs+1}$ to $x_{(r+1)s}$. By Observation (5), we can also conclude that $DF_r(2s)$ that we obtained admits a local antimagic 3-coloring with degree 2 (respectively, $3$ and $3s$) vertices having induced vertex labels $10k+1$ (respectively $13k+1$ and $s(17k+2)$). This completes the proof. \end{proof}

\begin{example} Using the twelve copies of $FB(1)$ as in Example~\ref{eg-FB12}, we can take $r=1$, $s=4$ to get $DF_1(8)$.  Note that the $DF_1(8)$ can be obtained from the top two figures in Example~\ref{eg-3DF4} by merging $y_1$ and $y_2$, and also $z_1$ and $z_2$. The $FB(4)$ component can be obtained from Example~\ref{eg-FB12} by merging vertices $x_i$, $i=5,6,7,8$. \rsq
\end{example}

\nt  For $1\le a\le s$, $1\le j\le r$, consider the $rDF(2s)$ as in Theorem~\ref{thm-rDF2seven}. Recall that the $j$-th component has
\begin{enumerate}[1.]\item $4s$ vertices of degree 2, namely $u_{(j-1)s+a}$, $u_{(2r-j)s+a}$, $v_{(j-1)s+a}$ and $v_{(2r-j)s+a}$, $1\le a\le s$,
\item  $2s$ vertices of degree 3, namely $w_{(j-1)s+a}$ and, $w_{(2r-j)s+a}$, $1\le a\le s$ and
\item 2 vertices of degree $3s$, namely $y_j$ and $z_j$.
\end{enumerate}

We define the following graphs.
\begin{enumerate}[1.]

\item $DF^1(r, 2s)$ is the graph obtained from $rDF(2s)$ by merging the degree 3 vertices in\break  $\{w_{(j-1)s+a}\;|1\le j\le r\}$ as $\alpha^{1,a}$ and in $\{w_{(2r-j)s+a}\;|\;1\le j\le r\}$ as $\alpha^{2,a}$.  Note that the degree of $\alpha^{i,a}$ is $3r$, for $1\le i\le 2$ and $1\le a\le s$.

\item $DF^2(r,2s)$ is the graph obtained from $rDF(2s)$ by merging the degree 2 vertices in\break  $\{u_{(j-1)s+a}\;|\; 1\le j\le r\}$ as $\beta^{1,a}$, in $\{u_{(2r-j)s+a}\;|\;1\le j\le r\}$ as $\beta^{2,a}$, in $\{v_{(j-1)s+a}\;|\; 1\le j\le r\}$ as $\beta^{3,a}$ and in $\{v_{(2r-j)s+a}\;|\;1\le j\le r\}$ as $\beta^{4,a}$. Note that the degree of $\beta^{i,a}$ is $2r$, for $1\le i\le 4$ and $1\le a\le s$.

\item $DF^3(r, 2s)$ is the graph obtained from $rDF(2s)$ by merging the degree $3s$ vertices in\break  $\{y_j\;|\;1\le j\le r\}$ as $y$ and in $\{z_j\;|\;1\le j\le r\}$ as $z$. Note that the degree of $y$ (and $z$) is $3rs$.
\end{enumerate}
Further, for even $r\ge 2$, denote by $DF^4(r,2s)$ the graph obtained from $rDF(2s)$ by merging $y_j$ to $z_{j+1}$ to obtain a degree $6s$ vertex, for $1\le j\le r$, where $z_{r+1}=z_1$. Note that $DF^3(2,2s) = DF^4(2,2s)$.

\begin{theorem}\label{thm-DF12(2s)} For $r\ge 2$, $i=1,2,3,4$, and $s\ge 1$,
\begin{enumerate}[1.]
\item $\chi_{la}(DF^1(r,2s)) = 3$ if $s$ is even and $rs\not\equiv 0\pmod 4$;
\item $\chi_{la}(DF^2(r,2s)) = 3$ if $s$ is even and $r\not\equiv 0\pmod 4$;
\item $\chi_{la}(DF^3(r,2s)) = 3$;
\item $\chi_{la}(DF^4(r,2s)) = 3$.
\end{enumerate}
 \end{theorem}

\begin{proof} Note that $DF^i(r,2s)$ is bipartite graph with equal partite sets size for each $i$, $1\le i\le 4$. By Lemma~\ref{lem-2part}, $\chi_{la}(DF^i(r,2s))\ge 3$. Suffice to show that $\chi_{la}(DF^i(r,2s))\le 3$ for $1\le i\le 4$.

\ms\nt Let $k=rs$. Using the $r$ copies of $DF(2s)$ in Theorem~\ref{thm-rDF2seven} and the corresponding local antimagic 3-coloring by merging the degree 3 vertices as defined for $DF^1(r,2s)$, we can immediately conclude that $DF^1(r,2s)$ admits a local antimagic 3-coloring having $2s$ (respectively, $4rs$ and $2r$) vertices of degree $3r$ (respectively, 2 and $3s$) with induced vertex label $r(13k+1)$ (respectively, $10k+1$ and $s(17k+2)$). From the proof of Theorem~\ref{thm-FBrs12}, we have shown that $r(13k+1)\ne s(17k+2)$ under the assumption.

\ms\nt Merging the degree 2 vertices as defined for $DF^2(r,2s)$, we can immediately conclude that $DF^2(r,2s)$ admits a local antimagic 3-coloring having $4s$ (respectively, $2rs$ and $2r$) vertices of degree $2r$ (respectively, 3 and $3s$) with induced vertex label $r(10k+1)$ (respectively, $13k+1$ and $s(17k+2)$). From the proof of Theorem~\ref{thm-FBrs12}, we have shown that $r(10k+1)\ne s(17k+2)$ under the assumption.

\ms\nt Merging the degree $3s$ vertices as defined for $DF^3(r,2s)$, we can immediately conclude that $DF^3(r,2s)$ admits a local antimagic 3-coloring having $4rs$ (respectively, $2rs$ and 2) vertices of degree 2 (respectively, 3 and $3rs$) with induced vertex label $10k+1$ (respectively, $13k+1$, and $rs(17k+2)$).

\ms\nt Merging the degree $3s$ vertices pairwise as defined for $DF^4(r,2s)$, we can immediately conclude that $DF^4(r,2s)$ admits a local antimagic 3-coloring having $r$ (respectively, $4rs$ and $2rs$) vertices of degree $6s$ (respectively, 2 and 3) with induced vertex label $2s(17k+2)$ (respectively, $10k+1$ and $13k+1$).

\ms\nt Thus, $\chi_{la}(DF^i(r,2s))\le 3$ for $i=1,2,3,4$.
\end{proof}

\section{$6\times 4n$ matrix}

We first construct a $6\times 4n, n\ge 1$ matrix using integers in $[1,20n]$ with integers in $[2n+1,4n]\cup [16n+1,18n]$ are used twice.
\[\fontsize{9}{11}\selectfont\begin{tabu}{|c|[1pt]c|c|c|c|c|[1pt]c|c|c|c|c|l}\cline{1-11}
R_1 &   1 &   2    & \cdots & 2n-1 & 2n &  6n+2 & 6n+4  & \cdots & 10-2 & 10 &\\\cline{1-11}
R_2 &16n+1 & 16n+2 & \cdots & 18n-1 & 18n & 18n &  18n-1  & \cdots & 16n+2 & 16n+1&\leftarrow\\\cline{1-11}
R_3 & 14n-1 & 14n-3 & \cdots & 10n+3 & 10n+1 & 6n &  6n-1  & \cdots & 4n+2 & 4n+1 &\\\cline{1-11}\cline{1-11}
R_4 & 14n+1 & 14n+2 & \cdots & 16n-1 & 16n & 6n+1 & 6n+3 & \cdots & 10n-3 & 10n-1 & \\\cline{1-11}
R_5 & 2n+1   & 2n+2 & \cdots & 4n-1 & 4n &  4n &   4n-1 & \cdots & 2n+2 & 2n+1 &\leftarrow\\\cline{1-11}
R_6 & 14n & 14n-2 & \cdots & 10n+4 & 10n+2 & 20n & 20n-1 & \cdots & 18n+2 & 18n+1 &\\\cline{1-11}
\end{tabu}\]
The `$\leftarrow$' indicates numbers in that row appear twice. We are now ready to trace $2n$ sequences, denoted $T_1, T_2, \ldots, T_{2n}$, of length 12 as follows:
\[\fontsize{6}{10}\selectfont\begin{tabu}{|c|[1pt]c|c|c|c|c|c|c|c|c|c|c|c|}\hline
T_1 & 1 & 16n+1 & 14n-1 & 6n+2 & 18n& 6n& 14n+1& 2n+1& 14n& 6n+1& 4n& 20n\\
T_2 & 2& 16n+2& 14n-3& 6n+4& 18n-1& 6n-1& 14n+2& 2n+2& 14n-2& 6n+3& 4n-1& 20n-1\\
\vdots & \vdots & \vdots & \vdots & \vdots & \vdots & \vdots & \vdots & \vdots & \vdots & \vdots & \vdots & \vdots \\
T_{n} & n & 17n & 12n+1& 8n& 17n+1& 5n+1& 15n& 3n& 12n+2& 8n-1& 3n+1& 19n+1\\\hline
T_{n+1} & 4n+1& 16n+1& 10n& 10n+1& 18n& 2n& 18n+1& 2n+1& 10n-1& 10n+2& 4n& 16n\\
T_{n+2}  & 4n+2& 16n+2& 10n-2& 10n+3& 18n-1& 2n-1& 18n+2& 2n+2& 10n-3& 10n+4& 4n-1 & 16n-1\\
\vdots & \vdots & \vdots & \vdots & \vdots & \vdots & \vdots & \vdots & \vdots & \vdots & \vdots & \vdots & \vdots \\
T_{2n} & 5n& 17n & 8n+2& 12n-1& 17n+1& n+1& 19n & 3n & 8n+1& 12n& 3n+1 & 15n+1\\\hline
\mbox{common} & +1 & +1 & -2 & +2 & -1 & -1 & +1 & +1 & -2 & +2 & -1 & -1 \\
\mbox{difference} & & \uparrow & & & \uparrow & & & \uparrow & & & \uparrow &
\end{tabu}\]
The `$\uparrow $' indicates numbers in that column appear twice. It is easy to verify that there is a bijective mapping between the entries in the matrix above and the terms of the sequences.

\ms\nt Consider the plane graph $C_{p}\times P_2$, where $p\ge 2$. We now define graph(s) $C_q(p,2)$ obtained from $C_{p}\times P_2$ by deleting $q$ edges not in $C_{p}$ and each pair of these $q$ edges are not incident to the same face, where $1\le q\le \frac{p}{2}$. Note that $C_q(p,2)$ is unique up to isomorphism if and only if $p=2q$ or $p=2q+1$ or $q=1$.

\begin{theorem}\label{thm-nC4(8,2)} For $s\ge 1$, $\chi_{la}(nC_4(8,2)) = 3$. \end{theorem}

\begin{proof} Let $G = nC_4(8,2)$. By definition, $G$ is a bipartite graph with equal partite set size. By Lemma~\ref{lem-2part}, $\chi_{la}(G)\ge 3$. Suffice to show that $G$ admits a local antimagic 3-coloring.

\ms\nt We now have the following observations.
\begin{enumerate}[(1)]
  \item For each $T_a$, $1\le a\le 2n$, the sum of the first and last terms (respectively, the $6$-th and 7-th; and the 9-th and 10-th terms) is a constant $20n+1$.
  \item For each $T_a$, $1\le a\le n$, the sum of the first 3 terms (respectively, the last 3 terms); whereas for $n+1 \le a\le 2n$, the sum of the 4-th to 6-th terms (respectively, the 7-th to 9-th terms), is a constant $30n+1$.
  \item For each $T_a$, $1 \le a\le n$, the sum of the 4-th to 6-th terms (respectively, the 7-th to 9-th terms); whereas for $n+1 \le a\le 2n$, the sum of the first 3 terms (respectively, the last 3 terms), is a constant $30n+2$.
  \item For $1\le a\le n$, $T_a$ and $T_{n+a}$ have the same second, fifth, eighth and eleventh terms.
\end{enumerate}

\nt For $1\le a\le n$, let the vertex set and the edge set of the $a$-th copy of $C_4(8,2)$ be\break $\{u_{a,i},v_{a,i}\;|\; 1\le i\le 8\}$ and $\{u_{a,i}u_{a,i+1}, v_{a,i}v_{a,i+1}, u_{a,2j}v_{a,2j}\;|\; 1\le i\le 8, 1\le j\le 4\}$, respectively, where $u_{a,9}=u_{a,1}$, $v_{a,9}=v_{a,1}$. We now define a bijective function $f : E(G) \to [1,20n]$ such that for the $a$-th copy of $C_4(8,2)$, $f$ assigns the edges of the $8$-cycle $u_{a,1}u_{a,2}\cdots u_{a,8}u_{a,1}$ (respectively, the $8$-cycle $v_{a,1}v_{a,2}\cdots v_{a,8}v_{a,1}$) by the first, third, fourth, sixth, seventh, ninth, tenth and twelfth terms of $T_a$ (respectively, $T_{n+a}$) consecutively, while the edges $u_{a,2}v_{a,2}$, $u_{a,4}v_{a,4}$, $u_{a,6}v_{a,6}$ and $u_{a,8}v_{a,8}$ are assigned by the second, fifth, eighth and eleventh terms of $T_a$. By the observations above, we can immediately conclude that $f$ is a local antimagic labeling of $G$ such that every degree 2 vertex has induced vertex label $20n+1$, where every two adjacent degree 3 vertices have induced vertex $30n+1$ and $30n+2$ respectively. Thus, $G$ admits a local antimagic 3-coloring. This completes the proof.
\end{proof}

\begin{example}\label{eg-6C4(82)} Take $n=6$, we have the sequences

{\fontsize{9}{12}\selectfont
\begin{multicols}{2}
\nt \begin{tabular}{>{$}r<{$}@{\ }l}
T_1: & 1, 97, 83, 38, 108, 36, 85, 13, 84, 37, 24, 120 \\
T_2: & 2, 98, 81, 40, 107, 35, 86, 14, 82, 39, 23, 119 \\
T_3: & 3, 99, 79, 42, 106, 34, 87, 15, 80, 41, 22, 118 \\
T_4: & 4, 100, 77, 44, 105, 33, 88, 16, 78, 43, 21, 117 \\
T_5: & 5, 101, 75, 46, 104, 32, 89, 17, 76, 45, 20, 116 \\
T_6: & 6, 102, 73, 48, 103, 31, 90, 18, 74, 47, 19, 115
\end{tabular}
\begin{tabular}{>{$}r<{$}@{\ }l}
T_7: &    25, 97, 60, 61, 108, 12, 109, 13, 59, 62, 24, 96 \\
T_8: &    26, 98, 58, 63, 107, 11, 110, 14, 57, 64, 23, 95 \\
T_9: &    27, 99, 56, 65, 106, 10, 111, 15, 55, 66, 22, 94 \\
T_{10}: & 28, 100, 54, 67, 105, 9, 112, 16, 53, 68, 21, 93 \\
T_{11}: & 29, 101, 52, 69, 104, 8, 113, 17, 51, 70, 20, 92 \\
T_{12}: & 30, 102, 50, 71, 103, 7, 114, 18, 49, 72, 19, 91
\end{tabular}
\end{multicols}}

\nt Then $6C_4(8,2)$ is labeled as follows:
\begin{figure}[H]
\centerline{\epsfig{file=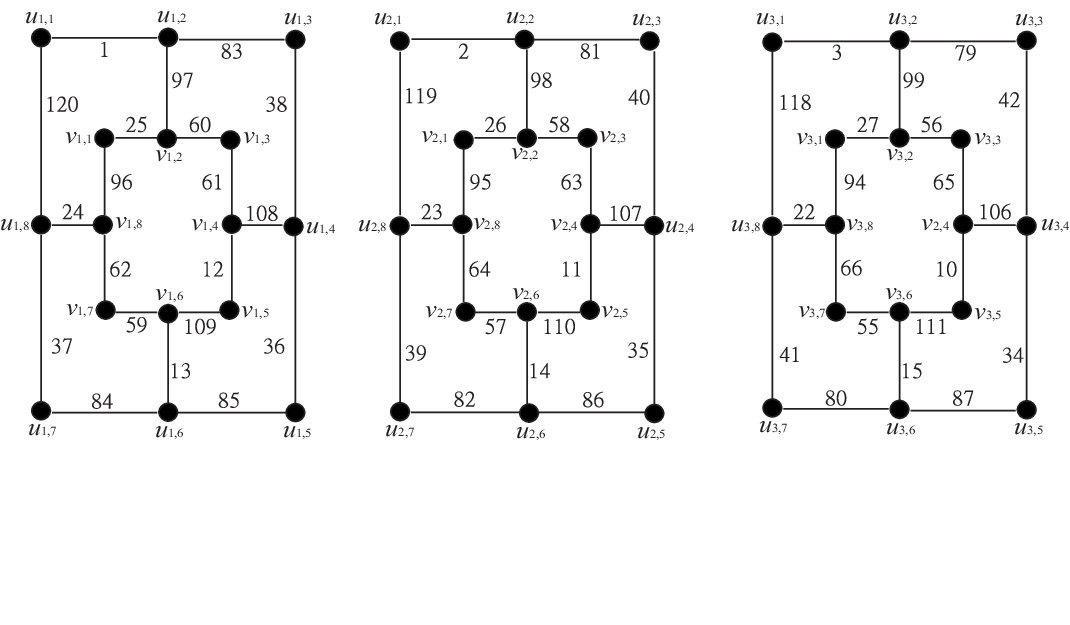, width=11.5cm}}
\vskip-1.8cm
\centerline{\epsfig{file=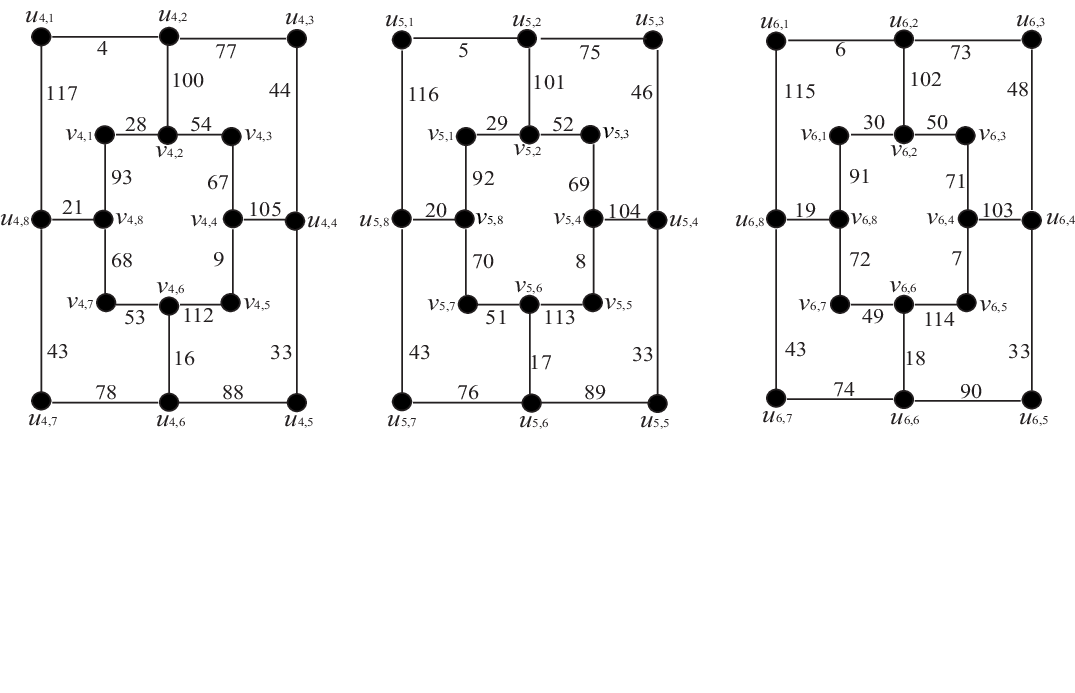, width=11.5cm}}
\vskip-3cm
\caption{The $6C_4(8,2)$ with the defined edge labeling.}
\end{figure}
\rsq
\end{example}

\nt We keep the names of vertices of $nC_4(8,2)$ defined in the proof of Theorem~\ref{thm-nC4(8,2)}. Suppose $n=rs$, $r\ge 1$, $s\ge 2$. Consider each $b\in [1,r]$. Let $G_1(r,s)$ be obtained from $rs$ copies of $C_4(8,2)$ by merging degree 2 vertices in $\{u_{(b-1)s+i,2j-1}\;|\; 1\le i\le s\}$ to get degree $2s$ vertices; and in $\{v_{(b-1)s+i,2j-1}\;|\; 1\le i\le s\}$ to get degree $2s$ vertices, for each $j\in[1,4]$, respectively.  Let $G_2(r,s)$ be obtained from $rs$ copies of $C_4(8,2)$ by merging degree 3 vertices in $\{u_{(b-1)s+i,2j} \;|\; 1\le i\le s\}$, for each $j=1,4$ respectively, to get degree $3s$ vertices; and in $\{v_{(b-1)s+i,2j}\;|\; 1\le i\le s\}$, for each $j=2,3$ respectively, to get degree $3s$ vertices. Note that $G_1(r,s)$, $G_2(r,s)$ are bipartite graphs with equal partite set size having $r$ component(s).

\begin{theorem}\label{thm-Gmrs} Suppose $m= 1,2$, $n=rs$, $r\ge 1$, $s\ge 2$, $\chi_{la}(G_m(r,s)) = 3$. \end{theorem}

\begin{proof} For $m=1,2$, since $G_m(r,s)$ is bipartite graph with equal partite set size, we have $\chi_{la}(G_m(r,s))\ge 3$. Suffice to show that $G_m(r,s)$ admits a local antimagic 3-coloring. Since $n=rs$, begin with the $nC_4(8,2)$ and the corresponding local antimagic 3-coloring as in the proof of Theorem~\ref{thm-nC4(8,2)}. Note that $G_1(r,s)$ has vertices of degree $2s$ and 3. By way of construction of $G_1(r,s)$, we immediately have a local antimagic labeling with induced vertex label of degree $2s$ vertices is $s(20n+1)$ and every two adjacent degree 3 vertices still have induced vertex labels $30n+1$ and $30n+2$ respectively.

\ms\nt Note that $G_2(r,s)$ has vertices of degree 2, 3 and $3s$.  Similarly, by way of construction of $G_2(r,s)$, we immediately have a local antimagic labeling with induced vertex label of degree 2 vertices is $20n+1$, degree 3 vertices is $30n+2$, and degree $3s$ vertices is $s(30n+1)$. Thus, $G_m(r,s)$, $m=1,2$ admits a local antimagic 3-coloring. This completes the proof.
\end{proof}

\begin{example} Using the first two components of $6C_4(8,2)$, we now give a component of the graphs $G_m(3,2)$ for $m=1,2$ as follows. Two other components can be obtained similarly.

\begin{figure}[H]
\centerline{\epsfig{file=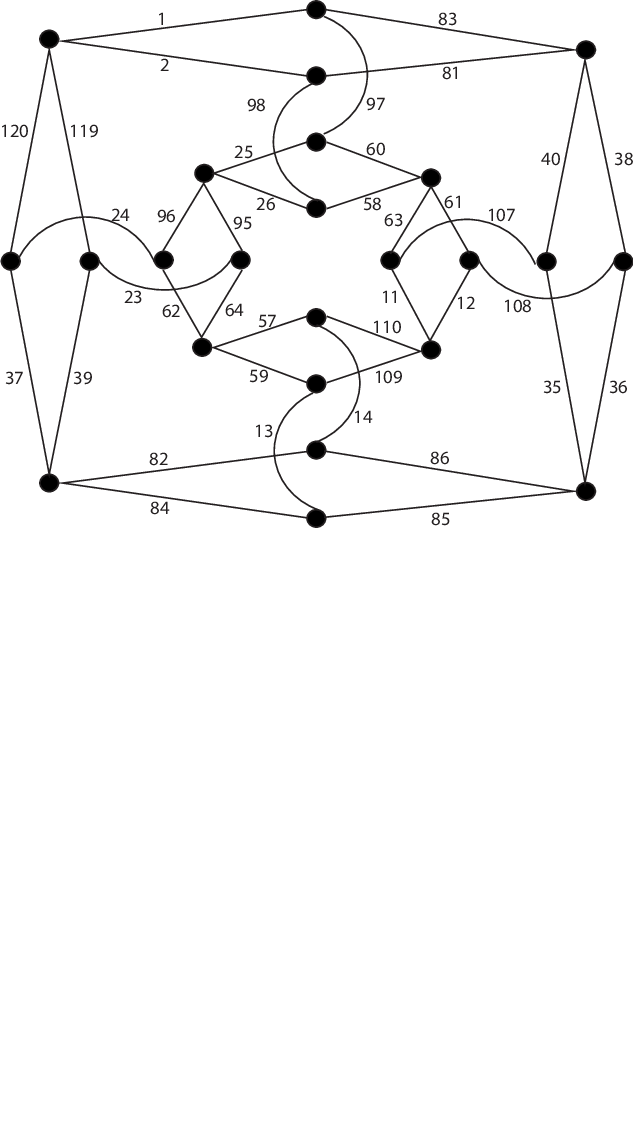, width=7cm}\hskip1cm\epsfig{file=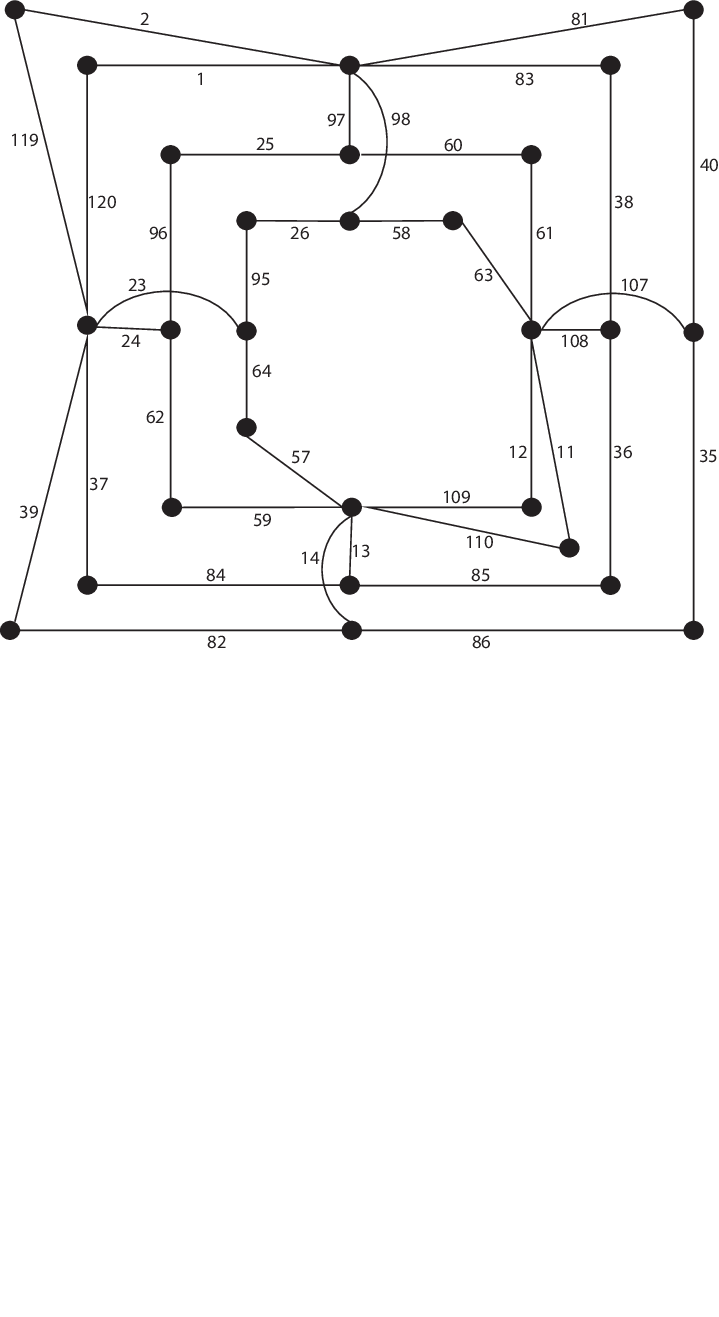, width=7cm}}
\vskip-6.cm
\caption{The first component of $G_m(3,2)$, $m=1,2$.}
\end{figure}
\rsq
\end{example}

\begin{enumerate}[(1)]
\item Let $H_1(n)$ be obtained from $n$ copies of $C_4(8,2)$ by merging degree 2 vertices in $\{u_{i,1}, u_{i,5}\}$ to get degree 4 vertex $x_{i,1}$; in $\{u_{i,3}, u_{i,7}\}$ to get degree 4 vertex $x_{i,2}$; in $\{v_{i,1}, v_{i,5}\}$ to get degree 4 vertex $y_{i,1}$; and in $\{v_{i,3}, v_{i,7}\}$ to get degree 4 vertex $y_{i,2}$, for each $i\in [1,n]$.
\item Let $H_2(n)$ be obtained from $n$ copies of $C_4(8,2)$ by merging the vertices in $\{u_{i,1}, v_{i,7}\}$ to get degree 4 vertex $x_{i,1}$; in $\{u_{i,5}, v_{i,3}\}$ to get degree 4 vertex $x_{i,2}$; in $\{u_{i,3}, v_{i,5}\}$ to get degree 4 vertex $y_{i,1}$; in $\{u_{i,7}, v_{i,1}\}$ to get degree 4 vertex $y_{i,2}$, for each $i \in [1,n]$.
\item Let $H_3(n)$ be obtained from $n$ copies of $C_4(8,2)$ by merging the vertices in $\{u_{i,1},v_{i,1}\}$ to get degree 4 vertex $x_{i,1}$; in $\{u_{i,5}, v_{i,5}\}$ to get degree 4 vertex $x_{i,2}$; in $\{u_{i,3},v_{i,3}\}$ to get degree 4 vertex $y_{i,1}$; in $\{u_{i,7},v_{i,7}\}$ to get degree 4 vertex $y_{i,2}$, for each $i\in [1,n]$.
\end{enumerate}

\nt Note that $H_1(n)$ is a bipartite graph with equal partite set size while $H_2(n)$ and $H_3(n)$ are tripartite graphs.

\begin{theorem}\label{thm-Hmn} Suppose $n\ge 1$, $m=1, 2, 3$, $\chi_{la}(H_m(n))=3$. \end{theorem}

\begin{proof} Since $H_1(n)$ is a bipartite graph with equal partite set size, we have $\chi_{la}(H_1(n))\ge 3$. Suffice to show that $H_1(n)$  admits a local antimagic 3-coloring.  Begin with the $nC_4(8,2)$ and the corresponding local antimagic 3-coloring as in the proof of Theorem~\ref{thm-nC4(8,2)}. By way of construction of $H_1(n)$, we immediately have a local antimagic labeling with induced vertex label of every two adjacent degree 3 vertices are $30n+1$ and $30n+2$ respectively, and the degree 4 vertices is $40n+2$.

\ms\nt Since $H_2(n)$ is tripartite, $\chi_{la}(H_2(n))\ge 3$. Similarly, by way of construction of $H_2(n)$, we immediately have a local antimagic labeling with induced vertex label of every degree 4 vertex is $40n+2$, and every two adjacent degree 3 vertices are $30n+1$ and $30n+2$ respectively. Thus, $H_m(n)$ $(m=1,2)$ admits a local antimagic 3-coloring.

\ms\nt Since $H_3(n)$ is tripartite, $\chi_{la}(H_3(n))\ge 3$. Similarly, by way of construction of $H_3(n)$, we immediately have a local antimagic labeling with induced vertex label of every degree 4 vertex is $40n+2$,and every two adjacent degree 3 vertices are $30n+1$ and $30n+2$ respectively. Thus, $H_m(n)$ $(m=1,2,3)$ admits a local antimagic 3-coloring.   This completes the proof.
\end{proof}


\begin{example} Using the first component of $6C_4(8,2)$, we now give the first component of the graphs $H_m(6)$ for $m=1,2,3$ as follows. Five other components can be obtained similarly.

\vskip-1.8cm
\begin{figure}[H]
\centerline{\epsfig{file=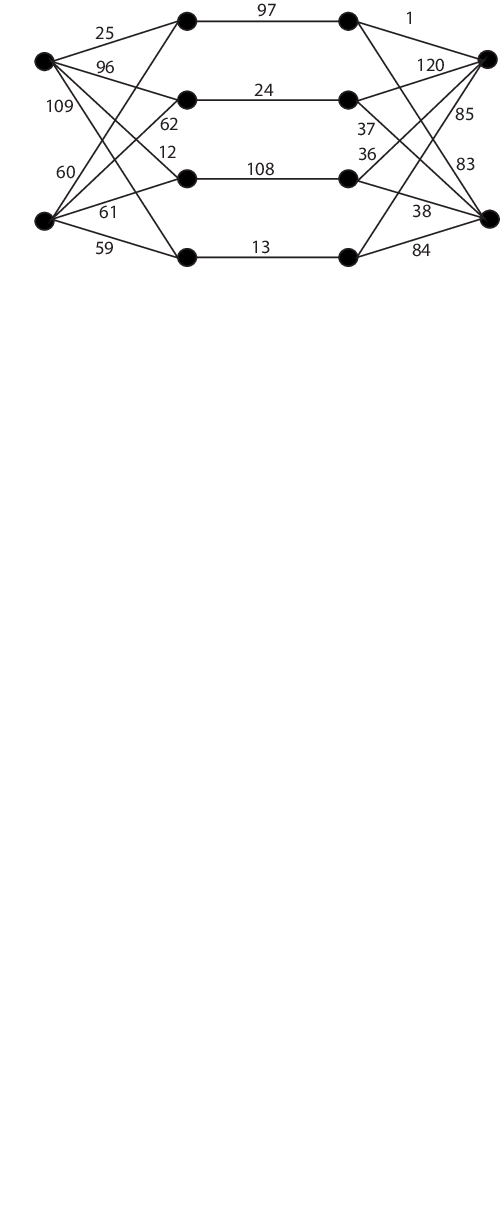, width=5cm}\hskip-4.5cm\epsfig{file=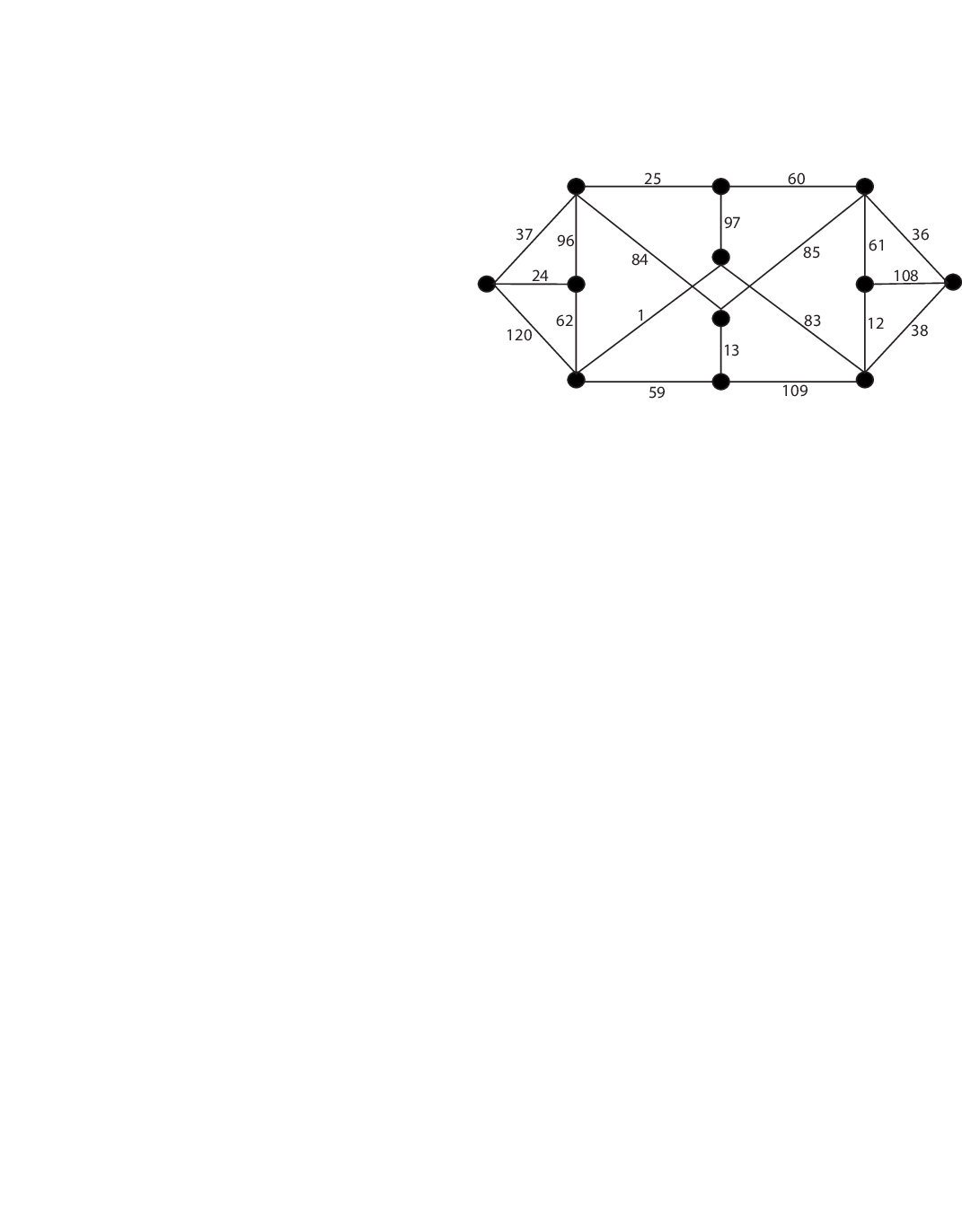, width=11cm}\hskip-4.5cm\vspace{1cm}\epsfig{file=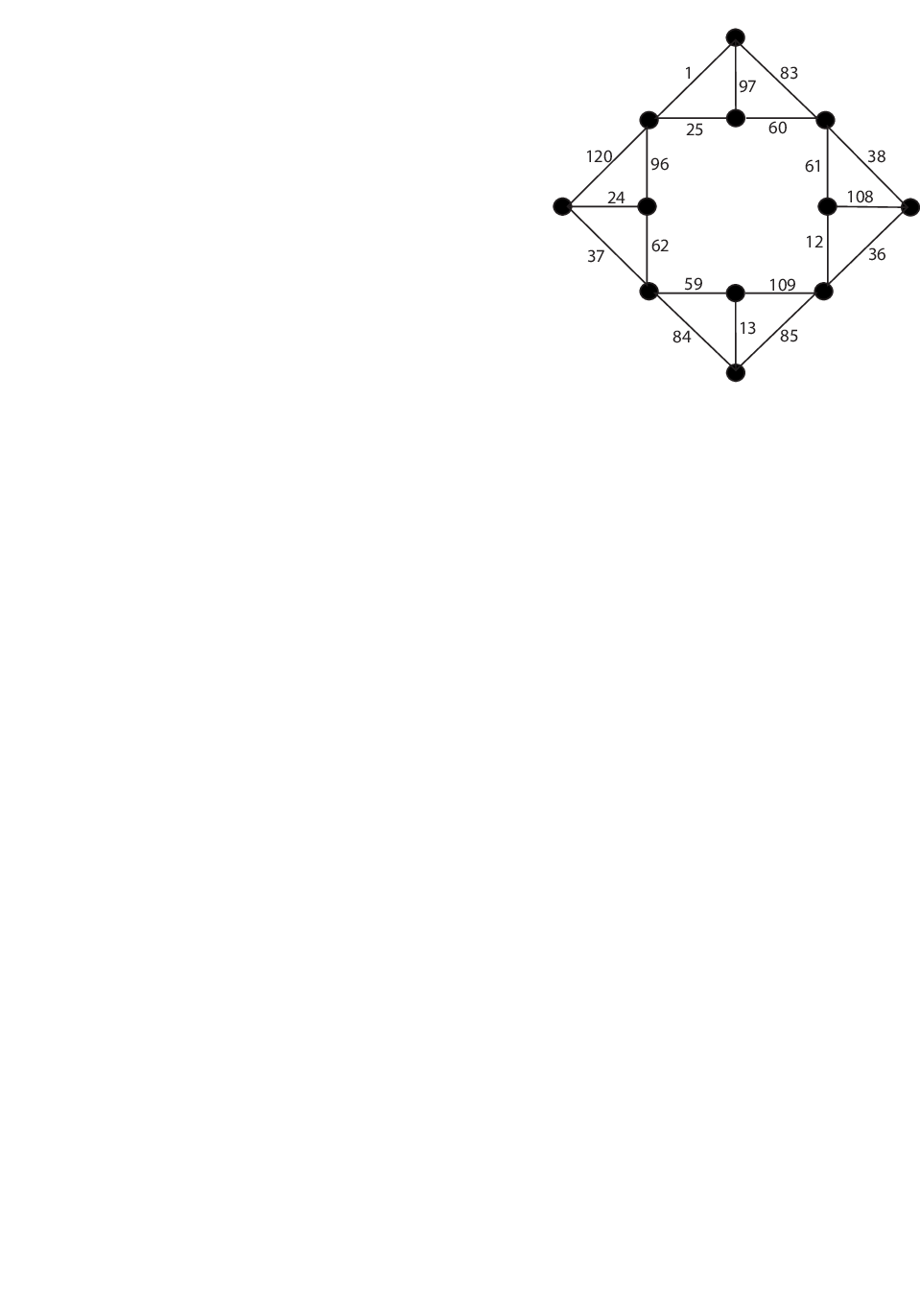, width=9cm}}
\vskip-10.cm
\caption{The first component of $H_m(6)$, $m=1,2,3$.}
\end{figure}
\rsq
\end{example}
%

\nt Suppose $n=rs$ with $r\ge 1$ and $s\ge 2$. For $m=1, 2, 3$  and each $b\in[1,r], j = 1,2$, let $H_m(r,s)$ be the graph obtained from $H_m(r,s)$ by merging degree 4 vertices in $\{x_{(b-1)s+i, j}\;|\;  i\in[1,s]\}$ to get degree $4s$ vertices; and in $\{y_{(b-1)s+i, j}\;|\;  i\in[1,s]\}$ to get degree $4s$ vertices.   Note that $H_1(r,s)$ is a bipartite graph with each component having equal partite set size so that $\chi_{la}(H_1(r,s)) \ge 3$. Moreover, $\chi_{la}(H_m(r,s))\ge \chi(H_m(r,s))=3$ for $m=2,3$. By using the same local antimagic $3$-coloring of $H_m(n)$, we get a new labeling for $H_m(r,s)$ with distinct induced vertex labels $s(40n+2)$, $30n+1$ and $30n+2$. So we have

\begin{theorem} For $r\ge 1$, $s\ge 2$ and $m=1,2,3$, $\chi_{la}(H_m(r,s))=3$.
\end{theorem}

\section{$k\times 10$ matrix}

\ms\nt We now construct a $k\times 10$ matrix for $k\ge 1$ (also with integers in $[1,10k]$) as follows to obtain all our remaining results.

\[\fontsize{7}{11}\selectfont
\begin{tabu}{|c|[1pt]c|c|c|c|c|c|c|c|[1pt]c|[1pt]c|}\hline
R_1 & 1 & 6k &  4k+1 & 2k+1 & 8k &  2k & 8k+1 & 10k & 7k & 3k+1 \\\hline
R_2 & k+1 & 6k+1  & 4k & 2k+2 & 8k-1 &  k  & 9k+1 & 9k & 6k-1& 4k+2 \\\hline
R_3 & k+2 & 6k+2 & 4k-1 & 2k+3 & 8k-2 &  k-1  & 9k+2 & 9k-1 & 6k-3 & 4k+4 \\\hline
\vdots  & \vdots & \vdots &\vdots  &\vdots  &\vdots & &\vdots&\vdots&\vdots&\vdots\\\hline
R_{k-1} & 2k-2 & 7k-2  & 3k+3 & 3k-1 & 7k+2 &  3 & 10k-2 & 8k+3 & 4k+5 & 6k-4 \\\hline
R_k  & 2k-1 & 7k-1 & 3k+2 & 3k & 7k+1 & 2 & 10k-1 & 8k+2  & 4k+3 & 6k-2\\\hline
\end{tabu}\]

\nt We now have the following observations:
\begin{enumerate}[(a)]
  \item In each row, sum of entries in columns $1, 8$ (respectively, $2,3$; $4, 5$; $6,7$ and $9,10$) is the constant $10k+1$.
  \item In each row, sum of entries in columns $1,2,9$ (respectively $5,6,10$) is the constant $13k+1$.
\end{enumerate}

\nt For $1\le i\le k$, $1\le j\le 8$, we use the row $i$ entries in columns 1 to 8 to construct the $i$-th copy of an 8-cycle with consecutive vertices $u_{i,1}, u_{i,2}, \ldots, u_{i,8}$ such that edge $u_{i,j}u_{i,j+1}$ is assigned the $j$-th entry, where $u_{i,9} = u_{i,1}$. We now add a new vertex $x_{i}$ adjacent to vertex $u_{i,2}$ and $u_{i,6}$. The columns 9 and 10 entries are now assigned to edges $x_iu_{i,2}$ and $x_iu_{i,6}$ respectively. Denote the graph obtained by $C^i_8$. Thus, for each $C^i_8$, $1\le i\le k$, vertex $u_{i,j}$, $j=1,3,5,7$ (respectively, $j=2,6$) are of degree 2 (respectively, 3) with induced vertex label $10k+1$ (respectively, $13k+2$). Moreover, vertices $x_i$, $u_{i,4}$ and $u_{i,8}$ are of degree 2 with induced vertex label $10k+1$, $6k+2$ and $18k+1$ respectively. Thus we have the following result:

\begin{theorem}\label{thm-Ck8} For $k\ge 1$, $2\le \chi_{la}(\sum\limits^k_{i=1} C^i_8) \le 4$. \end{theorem}

\begin{example}\label{eg-C4_82}  Take $k=4$, we get the following table and the $\sum\limits^4_{i=1} C^i_8$.
\[\fontsize{9}{11}\selectfont
\begin{tabu}{|c|[1pt]c|c|c|c|c|c|c|c|[1pt]c|[1pt]c|}\hline
R_1 & 1 & 24 &  17 & 9 & 32 &  8 & 33 & 40 & 28 & 13 \\\hline
R_2 & 5 & 25  & 16 & 10 & 31 &  4  & 37 & 36 & 23& 18 \\\hline
R_3 & 6 & 26 & 15 & 11 & 30 &  3  & 38 & 35 & 21 & 20 \\\hline
R_4 & 7 & 27 & 14 & 12 & 29 & 2 & 39 & 34  & 19 & 22\\\hline
\end{tabu}\]
\vspace*{-8mm}
\begin{figure}[H]
\centerline{\epsfig{file=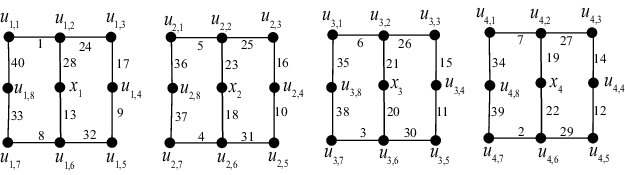, width=10cm}}
\caption{Graph $\sum^4\limits_{i=1} C^i_8$.}\label{fig:C4_82}
\end{figure}
\rsq
\end{example}

\nt Let $B^i_8$ be the bipartite graph obtained from $C^i_8$ by merging the vertices $u_{i,4}$ and $u_{i,8}$. Let $B_k = \sum\limits^k_{i=1} B^i_8$ for $k\ge 1$.

\begin{theorem} For $k\ge 1$, $2\le \chi_{la}(B_k)\le 3$. \end{theorem}

\begin{proof} Keeping the labeling as defined under Theorem~\ref{thm-Ck8}, we conclude that $B_k$ admits a local antimagic 3-coloring with induced vertex labels $10k+1, 13k+2, 24k+3$. This completes the proof. \end{proof}

\nt Let $C^i(8,2)$ be the graph obtained from $C^i_8$ by merging vertices $x_i$ and $u_{i,8}$. Let the new vertex be $z_i$. Thus, vertex $z_i$ is of degree 4 with induced vertex label $28k+2$. Since all $C^i(8,2)$ are isomorphic, we denote each of $C^i(8,2)$ by $C(8,2)$. 

\begin{theorem}\label{thm-kC82} For $k\ge 1$, $3\le \chi_{la}(kC(8,2)) \le 4$. \end{theorem}

\begin{proof} From  the observations above, we can immediately conclude that $kC(8,2)$ admits a local antimagic 4-coloring with each vertex $u_{i,j}, j=1,3,5,7$ of degree 2 (respectively, $u_{i,4}$ of degree 2, $u_{i,j}, j=2,6$ of degree 3 and $z_i$ of degree 4) has induced vertex label $10k+1$ (respectively, $6k+2$, $13k+1$ and $28k+2$). Thus, $\chi_{la}(kC(8,2))\le 4$. Since $\chi_{la}(kC(8,2))\ge \chi(kC(8,2)) = 3$, the theorem holds.  \end{proof}

\nt Let $D^i(8,2)$ be the graph obtained from $C^i_8$ by merging vertices $x_i$, $u_{i,4}$ and $u_{i,8}$. Let the new vertex be $w_i$. Thus, vertex $w_i$ is of degree 6 with induced vertex label $34k+2$. Since all $D^i(8,2)$ are isomorphic, we denote each of $D^i(8,2)$ by $D(8,2)$.

\begin{theorem} For $k\ge 1$, $\chi_{la}(kD(8,2))=3$.
\end{theorem}

\begin{proof} From the observations above, we can immediately conclude that $kD(8,2)$ admits a local antimagic $3$-coloring. Since $\chi(kD(8,2)= 3$, the theorem holds.
\end{proof}

\begin{theorem}\label{thm-evenrC82} Suppose $r\ge 1$ and $k = rs\ge 2$ for even $s\ge 2$. For $1\le i\le s/2$, $1\le a\le r$, if $G(8,2)$ is obtained from $s$ copies of $C^i(8,2)$ by merging vertices in $\{z_{(a-1)s+i}, u_{(2a-1)s/2+i,4}\}$ and in $\{z_{(2a-1)s/2+i}, u_{(a-1)s+i,4}\}$ respectively, then $\chi_{la}(rG(8,2))=3$.
\end{theorem}

\begin{proof} From the construction of $G(8,2)$, we know $\chi_{la}(rG(8,2))\ge \chi(rG(8,2)) = 3$. Suffice to show that $rG(8,2)$ admits a local antimagic 3-coloring. Observe that the edge labeling defined for $C^i_8$ $(1\le i\le k)$ corresponds to a local antimagic labeling of $rG(8,2)$ such that all the $4k$ vertices of degree 2 have induced vertex label $10k+1$, all the $2k$ vertices of degree 3 have induced vertex label $13k+1$, and all the $2r$ vertices of degree $3s$ have induced vertex label $s(17k+2)$. Thus, $\chi_{la}(rG(8,2))\le 3$. This completes the proof.  \end{proof}

\begin{example} From the $4C(8,2)$, we can take $r=1$ so that $G(8,2)$ is connected; or take $r=s=2$ so that $2G(8,2)$ has 2 components. The graph of the later case is shown below.
\begin{figure}[H]
\centerline{\epsfig{file=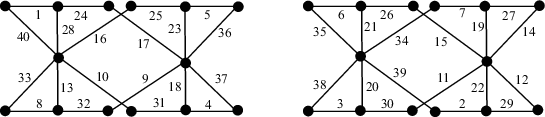, width=10cm}}
\caption{Graph  $2G(8,2)$.}\label{fig:2G82}
\end{figure} \rsq
\end{example}

\begin{theorem}\label{thm-oddkC82} Suppose $k=rs\ge 3$ is odd. For $1\le i\le (k+1)/2$ and $1\le j\le (k-1)/2$, let $H$ be obtained from $k$ copies of $C^i(8,2)$ by merging vertices in $\{z_i, u_{(k+1)/2+i}\}$ and in $\{z_{(k+1)/2+j}, u_{j,4}\}$, then $3\le \chi_{la}(H)\le 4$. \end{theorem}

\begin{proof} Clearly, $\chi_{la}(H)\ge \chi(H)=3$. 
%
%
Suffice to show that $\chi_{la}(H)\le 4$ for $k\ge 3$. From the construction of $H$, we observe that the edge labeling defined for $C^i_8$ $(1\le i\le k)$ corresponds to a local antimagic labeling of $H$ such that all $4k$ vertices of degree 2 have induced vertex label $10k+1$, all the $2k$ vertices of degree 3 have induced vertex label $13k+2$, and that there are $r$ vertices of degree $3(s-1)+4$ with induced vertex label $(s-1)(17k+2)+18k+1$ and another $r$ vertices of degree $3(s-1)+2$ with induced vertex label $(s-1)(17k+2)+6k+2$. Thus, $\chi_{la}(H)\le 4$. \end{proof}

\section{Conclusion and Open problems}

In this paper, we constructed three matrices of size $5\times 2k$, $6\times 2k$ ($k$ even) and $k\times 10$ respectively with entries in $[1,10k]$ that satisfy certain properties.  Consequently, many bipartite and tripartite graphs with local antimagic chromatic number 3 are obtained.  We have proved that $\chi_{la}(rFB(s)) = 3$ for odd $r,s\ge 3$ in~\cite[Theorem~2.2]{LSPN}.

\begin{conjecture} For even $r\ge 2$ and odd $s\ge 3$, $\chi_{la}(rFB(s))=3$.  \end{conjecture}

\ms\nt Note that each of the $n$ component of $G_3(n)$ in Theorem~\ref{thm-Hmn} is the triangular bracelet $TB(3)$ defined in~\cite[Theorem 3.2]{LSPN}. Thus, we have $\chi_{la}(k(TB(3))) = 3$ for $k\ge 1$. We end this paper with the following problems that arise naturally.

\begin{problem} Determine $\chi_{la}(k(TB(n))$ for $k\ge 2, n\ge 1$. \end{problem}

\begin{problem} For $k\ge 1$,  determine $\chi_{la}(B_k)$. \end{problem}

\begin{problem} For $k\ge 1$,  determine $\chi_{la}(\sum\limits^k_{i=1} C^i_8)$. \end{problem}

 \begin{problem} For $k\ge 1$, determine $\chi_{la}(kC(8,2))$. \end{problem}

\end{document}